\documentclass[12pt,psamsfonts,leqno,oneside,letterpaper]{amsart}
\usepackage[dvips,text={6.5truein,9truein},left=1truein,top=1truein]{geometry}
\usepackage{amssymb,amsmath,amscd,enumerate}
\usepackage[pdftex]{graphicx}
\usepackage{url}

\usepackage[colorlinks,linkcolor=blue,citecolor=blue,pdfstartview=FitH]{hyperref}
\input xy
\xyoption{all}
\SelectTips{cm}{12}

\usepackage{color}

\theoremstyle{definition}
\newtheorem{para}{}[section]

\newtheorem{remark}[para]{Remark}
\newtheorem{remarks}[para]{Remarks}
\newtheorem{notation}[para]{Notation}
\newtheorem{convention}[para]{Convention}
\newtheorem{definition}[para]{Definition}
\newtheorem{definitions}[para]{Definitions}

\newcommand\Alternatives{\begin{enumerate}[(i)]}
\newcommand\EndAlternatives{\end{enumerate}}
\newcommand\Conditions{\begin{enumerate}[(1)]}
\newcommand\EndConditions{\end{enumerate}}

\theoremstyle{plain}
\newtheorem{theorem}[para]{Theorem}
\newtheorem{lemma}[para]{Lemma}
\newtheorem{proposition}[para]{Proposition}
\newtheorem{corollary}[para]{Corollary}
\newtheorem{conjecture}[para]{Conjecture}
\newtheorem{claim}[equation]{}

\numberwithin{equation}{para}
\numberwithin{figure}{section}

\newcommand\Number{\begin{para}}
\newcommand\EndNumber{\end{para}}
\newcommand\Definition{\begin{definition}}
\newcommand\EndDefinition{\end{definition}}
\newcommand\Definitions{\begin{definitions}}
\newcommand\EndDefinitions{\end{definitions}}
\newcommand\Theorem{\begin{theorem}}
\newcommand\EndTheorem{\end{theorem}}
\newcommand\Conjecture{\begin{conjecture}}
\newcommand\EndConjecture{\end{conjecture}}
\newcommand\Remark{\begin{remark}}
\newcommand\EndRemark{\end{remark}}
\newcommand\Remarks{\begin{remarks}}
\newcommand\EndRemarks{\end{remarks}}
\newcommand\Convention{\begin{convention}}
\newcommand\EndConvention{\end{convention}}
\newcommand\Notation{\begin{notation}}
\newcommand\EndNotation{\end{notation}}
\newcommand\Lemma{\begin{lemma}}
\newcommand\EndLemma{\end{lemma}}
\newcommand\Proposition{\begin{proposition}}
\newcommand\EndProposition{\end{proposition}}
\newcommand\Corollary{\begin{corollary}}
\newcommand\EndCorollary{\end{corollary}}
\newcommand\Claim{\begin{claim}}
\newcommand\EndClaim{\end{claim}}
\newcommand\Proof{\begin{proof}}
\newcommand\EndProof{\end{proof}}
\newcommand\Equation{\begin{equation}}
\newcommand\EndEquation{\end{equation}}

\newcommand\Bullets{\begin{itemize}}
\newcommand\EndBullets{\end{itemize}}

\newcommand\op{{\mathfrak o}_{(\frakp)}}
\newcommand\opnought{{\mathfrak o}_{(\frakp_0)}}

\newcommand\trace{\mathop{\rm trace}}

\newcommand\RR{{\mathbb R}}

\newcommand\tGamma{\widetilde \Gamma}

\newcommand\tZ{\widetilde Z}

\newcommand\tQ{\widetilde Q}
\newcommand\tx{\widetilde x}

\newcommand\calc{{\mathcal C}}

\newcommand\calt{{\mathcal T}}

\newcommand\calg{{\mathcal G}}

\newcommand\calM{{\mathcal M}}
\newcommand\calo{{\mathcal O}}

\newcommand\caly{{\mathcal Y}}

\newcommand\hatGamma{\widehat\Gamma}
\newcommand\hatgamma{\widehat\gamma}

\newcommand\ZZ{{\mathbb Z}}

\newcommand\NN{{\mathbb N}}
\newcommand\CC{{\mathbb C}}
\newcommand\QQ{{\mathbb Q}}
\newcommand\HH{{\mathbb H}}

\newcommand\calf{{\mathcal F}}

\newcommand\frakp{{\mathfrak p}}
\newcommand\frakm{{\mathfrak m}}

\newcommand\Hom{\mathop{\rm Hom}}

\newcommand\vol{\mathop{\rm vol}}

\newcommand\pizzle{{\rm PSL}_2}
\newcommand\zzle{{\rm SL}_2}

\newcommand\ok{{\mathcal O}_K}

\newcommand\ooeese{\calo_{E,S^E}}
\newcommand\ooks{\calo_{K,S}}

\newcommand\ooksnought{\calo_{K_0,S_0}}

\begin{document}
\author{Peter B. Shalen}
\address{Department of Mathematics, Statistics, and Computer Science (M/C 249)\\  University of Illinois at Chicago\\
  851 S. Morgan St.\\
  Chicago, IL 60607-7045} \email{shalen@math.uic.edu}
\thanks{Partially supported by NSF grant DMS-0906155.}

\title{Orders of elements in finite quotients of Kleinian groups}

\begin{abstract}
A positive integer $m$ will be called a {\it finitistic order} for an element $\gamma$ of a group $\Gamma$ if there exist a finite group $G$ and a homomorphism $h:\Gamma\to G$ such that $h(\gamma)$ has order $m$ in $G$.
It is shown that up to conjugacy, all but finitely many elements of a given finitely generated, torsion-free Kleinian group admit a given integer $m>2$ as a finitistic order.
\end{abstract}

\maketitle

\section{Introduction}
I will be concerned with the following natural algebraic notion:

\Definition
Let $\gamma$ be an element of a group $\Gamma$. A positive integer $m$ will be called a {\it finitistic order} for $\gamma$ if there exist a finite group $G$ and a homomorphism $h:\Gamma\to G$ such that $h(\gamma)$ has order $m$ in $G$.
\EndDefinition

\Number\label{uneeda}
To illustrate this definition, consider the case in which $\Gamma$ is a free abelian group and $\gamma$ is a non-trivial element  of $\Gamma$. In this case there is an infinite cyclic direct summand $C$ of $\Gamma$ containing $\gamma$, and hence there is a homomorphism $h_0:\Gamma\to\ZZ$ such that $k:=h_0(\gamma)>0$. Given any positive integer $m$, the quotient homomorphism $\ZZ\to\ZZ/mk\ZZ$ maps $k$ onto an element of order $m$. Hence any positive integer $m$ is a finitistic order for $\gamma$ in this case.
\EndNumber

\Number
On the other hand, if $\gamma$ is an element of finite order $d$ in a group $\Gamma$, it is clear that only divisors of $d$ can be finitistic orders for $\gamma$. In particular, the only finitistic order for the identity element $1\in\Gamma$ is $1$. Likewise, if a group $\Gamma$ is not residually finite \cite{res-fin}, then by definition it contains at least one element $\gamma$ which is in the kernel of every homomorphism from $\Gamma$ to a finite group; the only finitistic order for such an element $\gamma$ is $1$.
\EndNumber

By a {\it Kleinian group} I will mean a discrete subgroup of $\pizzle(\CC)$. Such a group may be cocompact; it may be non-cocompact but have finite covolume; or it may have infinite covolume. The main result of this paper says that up to conjugacy, all but finitely many elements of a given finitely generated, torsion-free Kleinian group admit a given integer $m>2$ as a finitistic order. More precisely:

\Theorem\label{aathur}
Let $\Gamma$ be a finitely generated, torsion-free Kleinian group. Let $m>2$ be an integer, and let $X$ denote the set of all elements of $\Gamma$ for which $m$ is a finitistic order. Then $\Gamma-X$ is a union of finitely many conjugacy classes.
\EndTheorem

It is interesting to compare Theorem \ref{aathur} with some of the results proved in \cite{AllHam}. It follows from Lemmas 2 and 3 of \cite{AllHam}, together with  Proposition \ref{no great loss} of the present paper, that if $\Gamma$ is a finitely generated, torsion-free Kleinian group and $m$ is {\it any} positive integer, then for {\it every} element $\gamma$ of $\Gamma$ there is an integer {\it divisible} by $m$ which is a finitistic order for $\gamma$. Note that this result is neither stronger nor weaker than Theorem \ref{aathur}. 

The most novel ingredient in the proof of Theorem \ref{aathur} is
Proposition \ref{congker what} below, which is proved using a deep
number-theoretic result, due to Siegel and Mahler, about the
finiteness of the set of solutions to the $S$-unit equation in a
number field. Proposition \ref{congker what} implies that for an
arbitrary integer $m>2$ and an arbitrary finitely generated subgroup $\Gamma$ of
$\zzle(E)$, where $E$ is a number field, the {\it traces} of those
elements of $\Gamma$ that are not of finitistic order $m$ form a
finite set. The result is in fact stronger than this because the finite set of exceptional traces does not depend on the group $\Gamma$, but only the smallest number field $K$ containing the traces of all elements of $\Gamma$. 

Because Proposition \ref{congker what} does not require discreteness of the group $\Gamma$ and gives finiteness information based only the ``trace field'' $K$, it gives information that is not contained in Theorem \ref{aathur}. On the other hand, 
Proposition \ref{congker what} by itself does directly imply Theorem
\ref{aathur}, primarily because it establishes finiteness only for the set
of exceptional traces and not for  the set
of exceptional conjugacy classes, and secondarily because it requires $\Gamma$ to be contained in $\zzle(E)$ where $E$ is a number field, rather than in $\pizzle(\CC)$. Sections \ref{untitled section}---\ref{filling section} of this paper are devoted to the geometric arguments that are needed to deduce Theorem \ref{aathur} from the purely algebraic Propsoition \ref{congker what}. The various ingredients are assembled in Section \ref{gasoline alley} to prove Theorem \ref{aathur}.

Theorem \ref{aathur} is proved by establishing separate finiteness results for the loxodromic elements of a given finitely generated, torsion-free Kleinian group $\Gamma$ which are exceptional (in the sense that they do not admit a given integer $m>2$ as a finitistic order), and for the parabolic elements of $\Gamma$ which are exceptional. The finiteness for the exceptional loxodromic elements is reduced to Proposition \ref{congker what} via results about Kleinian groups established in Section \ref{GF section}. Among these results I would like to call attention to Proposition \ref{no great loss}, which I have not seen stated before in its general form; it is used to reduce the proof of finiteness for the exceptional loxodromic elements to the special case where the Kleinian group $\Gamma$ is geometrically finite, has no rank-$1$ maximal parabolic subgroups, and is contained in $\pizzle(E)$ for some number field $E$.

The work needed to prove finiteness for exceptional parabolic elements is done in Section \ref{filling section}. The general approach to the main result of this section, Proposition \ref{de uddah joik}, was directly inspired by the proof of \cite[Lemma 3]{AllHam}. The details of the proof turn out to involve Thurston's Dehn filling theorem and some interesting interactions between the topological and geometric aspects of hyperbolic $3$-manifolds.

Proposition \ref{congker what} will be applied in a different way to hyperbolic geometry in my forthcoming paper \cite{cubic}, which will establish interactions between the Margulis number of a hyperbolic $3$-manifold $M$ and the trace field of $M$.

I am grateful to Dick Canary and Steve Kerckhoff for some informative correspondence.

\section{Number fields, traces, and finitistic orders}

\Number\label{sweet perfection}
If $K$ is a field and $g$ is an element of $\zzle(K)$, I will denote by $[g]$ the image of $g$ under the quotient homomorphism $\zzle(K)\to\pizzle(K)$. If $Y$ is a subset of $K$ I will denote by $\calM_2(Y)$ the set of all $2\times2$ matrices whose entries lie in $Y$. If $Z$ is a subset of $\calM_2(K)$, I will denote by $\trace Z$ the set of all traces of elements of $Z$.
\EndNumber

\Lemma\label{newer dumb-diddly dumb}
Let $m>1$  be an integer, and $k$ be a finite field containing an element $\alpha$ whose order in the multiplicative group $k^\times$ is $2m$. Let $g$
be an element of $\zzle(k)$ such that $\trace g=\pm(\alpha+\alpha^{-1})$.  Then $[g]\in\pizzle(k)$ has order $m$.
\EndLemma

\Proof
First consider the case in which $\trace g=\alpha+\alpha^{-1}$. In this case, the characteristic polynomial $X^2-(\alpha+\alpha^{-1})X+1$ of $g$ has roots $\alpha$ and $\alpha^{-1}$ in $k$, and these roots are distinct since $\alpha$ has order $2m>1$. Hence $g$ is conjugate in $\zzle(k)$ to $\begin{pmatrix}\alpha&0\\0&\alpha^{-1}\end{pmatrix}$ and therefore has order $2m$ in $\zzle(k)$. Since $-I$ is the unique element of order $2$ in $\zzle(k)$ it follows that $[g]$ has order $m$ in $\pizzle(k)$.

In the case where $\trace g=-(\alpha+\alpha^{-1})$, we have $\trace (-g)=\alpha+\alpha^{-1}$. By the case already proved it follows that $[g]=[-g]$ has order $m$ in $\pizzle(k)$.
\EndProof

\Number\label{ooks defs}
In this
section I will use some concepts and elementary facts
from algebraic number theory that will be used.
The book \cite{neukirch} is a general reference.

Let $K$ be a number field. I shall denote by $S_\infty(K)$ the set of
all archimedean places of $K$. 

A set $S$ of places of $K$ will be termed {\it admissible} if $S$ is
finite and $S\supset S_\infty(K)$. If $S$ is admissible, I shall
denote by $\calo_{K,S}$ the {\it ring of $S$-integers} of $K$, defined
as the intersection of all valuation rings corresponding to places not
belonging to $S$. In particular, $\ok:=\calo_{K,\emptyset}$  is the
ring of integers of $K$ (see the statement ($\alpha$) on p. 264 of \cite{rib}).
\EndNumber

\Number\label{it don't matter}  If
$K$ and $E$ are number fields, with $K\subset E$, and if $S$ is an
admissible set of places in a number field $K$, I shall denote by
$S^E$ the set of all extensions to $E$ of places in $S$. Since a given
valuation of $K$ admits only finitely many extensions to valuations of
$E$, for example by \cite[Proposition 8.2]{neukirch}, $S^E$ is an
admissible set of places of $E$. On the other hand, since every
valuation of $K$ admits at least one extension to a valuation of
$E$---again for example by \cite[Proposition 8.2]{neukirch}---we have
$\ooks=\ooeese\cap K$.
\EndNumber

\Number\label{taragon and garlic}
Let $E$ be a number field,  and let $\frakp$ be a nonarchimedean place of $E$. I will
denote by $\op$ the valuation ring defined by $\frakp$, by  $k_\frakp$  the
residue field $\op/\frakp$, and by $\eta_\frakp:\op\to k_\frakp$ the quotient
homomorphism. I will denote by
$h_\frakp$ the natural homomorphism $\zzle(\op)\to\zzle(k_\frakp)$, defined by
$$\begin{pmatrix}a&b\\c&d\end{pmatrix}\mapsto\begin{pmatrix}\eta_\frakp(
a)&\eta_\frakp( b)\\\eta_\frakp( c)&\eta_\frakp( d)\end{pmatrix}.$$
I will denote by
$\kappa_\frakp:\zzle(\op)\to\pizzle(k_\frakp)$
the homomorphism defined by $\kappa_\frakp(\gamma )=[h_\frakp(\gamma )]$.
\EndNumber

\Lemma\label{dreamotomic}
For each integer $n>1$ there is a positive integer $N=N(n)$ with the following property.
Let $E$ be any number field, and let $\frakp$ be a nonarchimededean place of $E$ such that $1/N\in\op$. Let $\omega\in E$ be a root of unity of order $n$. Then $\omega\in \calo_{E}\subset\op$, and $\eta_\frakp(\omega)$ is an element of order $n$ in the multiplicative group $k_\frakp^\times$.
\EndLemma

\Proof
For each positive integer $m$  the cyclotomic polynomial $\Phi_m(X)\in\ZZ[X]$ is irreducible in $\QQ[X]$ and has degree $m$.
For each proper divisor $d$ of $n$, the polynomial $X^d-1$ is a product of cyclotomic polynomials whose degrees at most $d<n$. Hence $X^d-1$ is relatively prime to $\Phi_n(X)$ in $\QQ[X]$, and so there exist polynomials $A_d(X)$ and $B_d(X)$ in $\ZZ[X]$ such that $A_d(X)\Phi_n(X)+B_d(X)(X^d-1)=N_d$ for some $N_d\in\QQ$. I will take $N$ to be a positive integer which is divisible for $N_d$ for every proper divisor $d$ of $n$, and show that the conclusion holds with this choice of $N$.

Let $E$ be any number field, and let $\frakp$ be a nonarchimededean place of $E$ such that $1/N\in\op$. Let $\omega\in E$ be a root of unity of order $n$.
Then in particular $\omega$ is an algebraic integer, and so $\omega\in \calo_{E}\subset\op$. Let us set $\alpha=\eta_\frakp(\omega)$. Since $\eta_\frakp:\op\to k_\frakp$  is homomorphism, we have $\alpha^n=1$. I shall complete the proof by assuming that the order of $\alpha$ is a proper divisor $d$ of $n$ and deriving a contradiction.

Since the equality $A_d(X)\Phi_n(X)+B_d(X)(X^d-1)=N_d$ holds in
$\ZZ[X]$, and since $N_d|N$ and $1/N\in\op$, we have
$A_d(\alpha)\Phi_n(\alpha)+B_d(\alpha)(\alpha^d-1)\ne0$ in
$k_\frakp$. But since $\omega$ has order $n$ in $E^\times$, we have
$\Phi_n(\omega)=0$, and since $\eta_\frakp$ is a homomorphism we have
$\Phi_n(\alpha)=0$. Furthermore, since $\alpha$ has order $d$ in
$k_\frakp^\times$, we have $\alpha^d-1=0$. This is a
contradiction.  \EndProof

I will now turn to the main result of the section.

\Proposition\label{congker what}
Let $K$ be a number field, let $S$ be an admissible set of places of
$K$, and let $m>2$ be an integer. Then there is a finite set
$W\subset\ooks$ with the following property. Let $\gamma$ be an element of
$\zzle(\ooks)$ such that $\trace\gamma\notin W$. Then there exists a
place $\frakp$ of $K$, with $\frakp\notin S$, such that $\kappa_\frakp(\gamma)\in \pizzle(k_\frakp)$ has
order $m$.
 In particular, $m$ is a finitistic order for
$\gamma\in\zzle(\ooks)$. 
\EndProposition

\Proof 
Set $n=2m$. Let $N=N(n)$ be a positive integer having the property stated in Lemma \ref{dreamotomic}.
Let $K_0$ be a finite extension of $K$ which contains a primitive
$n$-th root of unity $\omega$, and set $\tau=\omega+\omega^{-1}\in
K_0$. Since $m>2$ we have $\omega^2\ne-1$ and hence $\tau\ne0$. Let
$S_0$ denote an admissible set of places of $K_0$, containing $S^{K_0}$
(see \ref{it don't matter}), such that $2$, $N$ and $\tau $ are units in $\ooksnought$. 

Let $U$ denote the set of all elements $u\in\ooksnought$ such that both $u$ and $1-u$ are units in $\ooksnought$. According to 
\cite[Theorem D.8.1]{hs} (a result due to Siegel and Mahler), $U$ is a
finite set. Hence the set $W_0:=\{(2u-1)\tau :u\in U\}\subset K_0$  is also finite, and so is $W=K\cap W_0$. I shall
complete the proof of the proposition by showing that the conclusion is true with this choice of $W$.

I claim that for any $t\in K_0-W_0$, either $t-\tau $ or $t+\tau $ is a non-unit in $\ooksnought$. Indeed, suppose that $t-\tau $ and $t+\tau $ are both units. Since $2$ and $\tau $ are also units in $\ooksnought$, it then follows that $u:=(\tau +t)/2\tau $ and $1-u=(\tau -t)/2\tau $ are units as well, i.e. that $u\in U$. Hence $t= (2u-1)\tau \in W_0$, a contradiction.

Now let  $\gamma$ be any element of
$\zzle(\ooks)$ such that $t:=\trace\gamma\notin W$. In particular we have $t\in K_0-W_0$. Hence
  we may choose an element $\epsilon$ of $\{1,-1\}$ such that
  $t-\epsilon\tau $ is a non-unit in $\ooksnought$.

Set $S_0=S^{K_0}$. By \ref{it don't matter} we have $\ooksnought\cap K=\ooks$, and hence $t-\epsilon\tau $  is a non-unit in $\ooksnought$. 
 It follows that there is a place  $\frakp_0$ of $K_0$ such that
$t-\epsilon\tau $ belongs to the maximal ideal $\frakm_0$   of $\opnought$. 
Hence, in the notation of \ref{taragon and garlic}, we have $\eta_{\frakp_0}(t)=\epsilon \eta_{\frakp_0}(\tau )\in k_{\frakp_0}$. But
since $t=\trace\gamma$ we have $\eta_{\frakp_0}(t)=\trace
h_{\frakp_0}(\gamma)$. Thus the element $g:=h_{\frakp_0}(\gamma)$ of $\zzle(k_{\frakp_0} )$ has trace $\epsilon \eta_{\frakp_0}(\tau)$. If we set 
$\alpha=\eta_{\frakp_0}(\omega)$, then $\eta_{\frakp_0}(\tau)=\alpha+\alpha^{-1}$. 
Since $\omega$ is a primitive $n$-th root of unity, and since $1/N\in\ooksnought\subset\opnought$, it follows from Lemma \ref{dreamotomic}  that $\alpha$ has order $n=2m$ in $k_{\frakp_0}^\times$. Since $\trace g=\epsilon( \alpha+\alpha^{-1})$, it follows from Lemma \ref{newer dumb-diddly dumb} that $[g]=\kappa_{\frakp_0}(\gamma)\in\pizzle(k_{\frakp_0})$ has order $m$. 

We have $\ker \kappa_{\frakm_0}=\{\pm I\}+\calM_2(\frakm_0)$. If
$\frakp$ denotes the restriction of the place $\frakp_0$ to $K$, and
$\frakm$ denotes the maximal ideal of $\op$, we have $\ker
\kappa_{\frakp}=\{\pm I\}+\calM_2(\frakm)$. Hence $\ker
\kappa_{\frakp}=\zzle(K_0)\cap \ker \kappa_{\frakp_0}$. It follows
$\kappa_{\frakp}(\gamma) $ has the same order as
$\kappa_{\frakp_0}(\gamma) $, namely $m$.
\EndProof

\section{Preliminaries on Kleinian groups}\label{untitled section}
The study of torsion-free Kleinian groups is closely related to the study of (complete) orientable hyperbolic $3$-manifolds. The topological group $\pizzle(\CC)$ may be identified by a continuous isomorphism with the group of orientation-preserving isometries of the $3$-dimensional hyperbolic space $\HH^3$. If $\Gamma\le\pizzle(\CC)$ is discrete and torsion-free, the action of $\Gamma$ on $\HH^3$ is free and properly discontinuous, and $M:=\HH^3/\Gamma$ inherits the structure of an orientable hyperbolic $3$-manifold with $\pi_1(M)\cong\Gamma$.
Up to conjugacy there is a natural identification of $\pi_1(M)$ with $\Gamma$.

In this section I will collect a few essential facts about hyperbolic $3$-manifolds which will be needed in Sections \ref{GF section} and \ref{filling section}.

\Number\label{general stuff reporting for duty}
I will follow the conventions of \cite[Section 3]{bounds}. In particular, I will work in the smooth category (so that manifolds and submanifolds are understood to be smooth) but will often quote results proved in the piecewise linear category; the justification for doing this is explained in \cite[Subsection 3.1]{bounds}. It is understood that a {\it manifold} may have a boundary. Recall from \cite[Section 3]{bounds} that a $3$-manifold $M$ is said to be {\it irreducible} if $M$ is connected and every $2$-sphere in $M$ is the boundary of a $3$-ball in $M$. According to \cite[Proposition 3.8]{bounds}, every orientable hyperbolic $3$-manifold is irreducible.

I will also follow the conventions of \cite[Section 2]{bounds} in statements and arguments involving fundamental groups: I will suppress base points
whenever it is possible to do so without ambiguity.

An element of $\pizzle(\CC)$ is said to be {\it parabolic} if it is non-trivial and has the form $[A]$ (see \ref{sweet perfection}) for some $A\in\zzle(\CC)$ with $\trace A=2$. By a {\it parabolic subgroup} of a Kleinian group  $\Gamma$ I will mean a non-trivial subgroup of  $\Gamma$ whose non-trivial elements are all parabolic.

By a {\it standard cusp neighborhood} $X$ in the orientable hyperbolic $3$-manifold $M=\HH^3/\Gamma$ I will mean a subset of the form $B/\Gamma_{B}$, where $B$ is an open horoball in $\HH^3$, precisely invariant under $\Gamma$ (in the sense that $\gamma\cdot B$ is either equal to $B$ or disjoint from $B$ for every $\gamma\in\Gamma$), and the stabilizer $\Gamma_B$ of $B$ is a parabolic subgroup. I will define the {\it rank} of $X$ to be the rank of $\Gamma_B$ (which must be equal to $1$ or $2$).

In general a Kleinian group is said to be {\it elementary} if it has an abelian subgroup of finite index. A torsion-free elementary Kleinian group is itself abelian (see, e.g., \cite[Proposition 2.1]{finiteness}).

\EndNumber

\Number \label{vexcor}
A non-elementary, orientable hyperbolic $3$-manifold $M$ has a well-defined {\it convex core}, which I will denote by $\calc(M)$ . By definition, $\calc(M)$ is the smallest non-empty closed subset of $M$ which is {\it convex} in the strong sense that every geodesic path with endpoints in $\calc(M)$ is entirely contained in $\calc(M)$. For a construction of $\calc(M)$ and proofs of its basic properties, see \cite[Subsection 
3.1.1]{M-T} or \cite[p. 63]{morgan}. I will denote by $\calc_1(M)$ the closed radius-$1$ metric neighborhood of $\calc_1(M)$; it follows from \cite[Proposition 3.1]{M-T}  that $\calc_1(M)$ is a $3$-manifold and a deformation retract of $M$. 

A {\it closed geodesic} in a hyperbolic manifold $M$ will be regarded as a map $C:S^1\to M$ such that the map $t\mapsto C(e^{2\pi it})$ from $\RR$ to $M$ is a geodesic; I will write $|C|=C(S^1)$. The construction of $\calc(M)$ given in \cite{M-T} or \cite{morgan} immediately implies that $|C|\subset\calc(M)$ for any closed geodesic $C$.
\EndNumber

\Number\label{re-verse it}
As in \cite[p. 55]{M-T}, I will define the $\epsilon${\it-thin part} $M_{(0,\epsilon)}$ of a non-elementary orientable hyperbolic $3$-manifold $M$, where $\epsilon>0$ is given, to be the set of all points of $M$ which are base points of homotopically non-trivial loops of length less than $\epsilon$. According to \cite[Theorem 2.24]{M-T}, there is a universal constant $\epsilon_0$ such that for every orientable $3$-manifold $M$, each component of $M_{(0,\epsilon_0)}$ is either a standard cusp neighborhood or a metric neighborhood of a simple closed geodesic. I will fix such a constant $\epsilon_0$ for the rest of the paper. According to \cite[Lemma 6.7]{morgan}, each rank-$2$ standard cusp neighborhood $X$ in $M$ contains a smaller  rank-$2$ standard cusp neighborhood $s(X)$ which is contained in $\calc(M)$. I will refer to the set $K:=\calc_1(M)-\bigcup_X s(X)$, where $X$ ranges over all 
components of $M_{(0,\epsilon_0)}$ which are rank-$2$ standard cusp neighborhoods, as a {\it truncation} of $\calc_1(M)$. Thus every non-elementary orientable hyperbolic $3$-manifold admits a truncation $K$. (According to this definition, a truncation of $\calc_1(M)$ is not quite unique, as it depends on the choice of the standard cusp neighborhoods $s(X)$.)
\EndNumber

\Number\label{erhard}
Let  $K$ be an orientable, irreducible $3$-manifold.
An {\it essential singular torus} in $K$ is defined to be a map $f:T^2\to K$ such that (i) $f_\sharp:\pi_1(T^2)\to\pi_1(K)$ is injective, and (ii) $f$ is not homotopic in $K$ to a map of $T^2$ into $\partial K$.

Now suppose that $\calt$ is a $2$-dimensional submanifold of $\partial K$ such that the inclusion homomorphism $\pi_1(T)\to\pi_1(K)$ is injective for every component $T$ of $\calt$. In this setting I will define an {\it essential singular annulus} in the pair $(K,\calt)$ to be a
map of pairs $f:(S^1\times[0,1],\partial(S^1\times[0,1]))\to (K,\calt)$ such that
$f_\sharp:\pi_1(S^1\times[0,1])\to\pi_1( K)$ is injective, and $f$ is
not homotopic rel $\partial(S^1\times[0,1])$ to a map of
$S^1\times[0,1])$ into $\calt$.)
 
\EndNumber

\Proposition\label{lame song}
Let $M=\HH^3/\Gamma$ be an orientable hyperbolic $3$-manifold, and let $K$ be a truncation of $\calc_1(M)$. Then:
\begin{enumerate}
\item $K$ is irreducible and has no essential singular tori, and is a deformation retract of $M$. 
\item For every torus component $T$ of $\partial K$ there is a standard cusp neighborhood which is a component of $\overline{M-K}$ and is bounded by $T$; furthermore, the inclusion homomorphism $\pi_1(T)\to\pi_1(M)$ is injective, and its image, a subgroup of $\pi_1(M)=\Gamma$ defined up to conjugacy, is a rank-$2$ maximal parabolic subgroup of $\Gamma$. 
\item Conversely, every rank-$2$ maximal parabolic subgroup of $\Gamma$ is conjugate to
the image of the inclusion homomorphism $\pi_1(T)\to\pi_1(M)$ for some torus component $T$  of $\partial K$.
\end{enumerate}
\EndProposition

\Proof 
Let $p:\HH^3\to M$ denote the quotient projection.

Suppose that $P$ is a rank-$2$ maximal parabolic subgroup of $\Gamma$. Fix generators $\gamma_1$ and $\gamma_2$ of $P$. Since $\gamma_1$ and $\gamma_2$ are commuting parabolics, there is a point $\tx\in\HH^3$ such that $d(\tx,\gamma_i\cdot\tx)<\epsilon_0$ for $i=1,2$. Set $x=p(\tx)$. The base point $\tx\in p^{-1}(x)\subset\HH^3$ determines an isomorphic identification of $\pi_1(M,x)$ with $\Gamma$, and under this identification, each $\gamma_i$ is represented by a loop $c_i:[0,1]\to M$ of length $<\epsilon_0$. If we set $G=c_1([0,1])\cup c_2([0,1])$, then each point of $G$ is the basepoint of a loop of length $<\epsilon_0$, and hence $G\subset M_{(0,\epsilon_0)}$. Under our identification, $P$ in contained in the image $P'$ of the inclusion homomorphism $\pi_1(G,x)\to\pi_1(M,x)$, and hence in the image of the inclusion homomorphism $\pi_1(X,x)\to\pi_1(M,x)$, where $X$ denotes the component of $M_{(0,\epsilon_0)}$ containing $G$. Since $P$ is non-cyclic, we the discussion in \ref{re-verse it} shows that $X$ is a rank-$2$ standard cusp neighborhood; hence $P'$ is a parabolic subgroup of $\Gamma$, and the maximality of $P$ implies that $P=P'$. It now follows from the definition of a truncation that $X$ contains a unique torus component $T$ of $\partial K$, and that the image (a priori defined up to conjugacy) of the inclusion homomorphism $\pi_1(T)\to\pi_1(K)$ is conjugate to $P$. This proves Assertion (3).

Conversely, suppose that $T$ is a torus component of $\partial K$. According to the definition of a truncation, $T$ is either (a) a component of $\partial\calc_1(M)$, or (b) the boundary of 
a standard cusp neighborhood which is a component of $\overline{M-K}$. However, by \cite[Proposition 3.1]{M-T}, $\partial\calc_1(M)$ is homeomorphic to $\Omega/\Gamma$, where $\Omega$ denotes the set of discontinuity of $\Gamma$. 
Hence if (a) holds, some component of $\Omega$ conformally covers a torus; this is impossible for a non-elementary Kleinian group $\Gamma$, because the limit set $\Lambda=\widehat\CC -\Omega$ must contain more than two points, and hence each component of $\Gamma$ is a hyperbolic Riemann surface. Hence (b) must hold. If $X$ denotes the standard cusp neighborhood bounded by $T$, then by definition $p^{-1}(X)$ is a horoball $B\subset\HH^3$, and $p^{-1}(T)$ is the frontier $H$ of $B$ in $\HH^3$. Since $H$ is simply connected, the inclusion homomorphism $i:\pi_1(T)\to\pi_1(M)$ is injective. Furthermore, the image of is identified with the stabilizer $\Gamma_B$ of $B$ in $\Gamma$, which is a maximal parabolic subgroup of $\Gamma=\pi_1(M)$; being isomorphic to $\pi_1(T)$, it has rank $2$. This proves (2).

To prove (1), first note that according to the definition of a truncation, each component of $\calc_1(M)-K$ is a standard cusp neighborhood bounded by a torus component of $\partial K$; hence $K$ is a deformation retract of $\calc_1(M)$. Since $calc_1(M)$ is in turn a deformation retract of $M$ (cf. \ref{vexcor}, it follows that $K$ is a deformation retract of $M$.

Next note that, by \cite[Proposition 3.1]{M-T},  $\calc_1(M)$ is a deformation retract of $M$; hence every component of $M-\calc_1(M)$ has non-compact closure in $M$. In view of the definition of a truncation, it follows that every component of $M-K$ has non-compact closure in $M$.  Since $M$ is irreducible by \cite[Proposition 3.8]{bounds}, it follows that $K$ is irreducible. Finally, suppose that $f:T^2\to K\subset M$ is an essential singular torus. Then $f_\sharp(\pi_1(T^2)\le\pi_1(M)$ is a rank-$2$ free abelian subgroup of $\pi_1(M)$ which is defined up to conjugacy. There is an isomorphic identification of $\pi_1(M)$ with $\Gamma$ which is also canonically defined up to conjugacy. Since $\Gamma$ is discrete, $X$ must be parabolic. Let $X_0$ be a maximal parabolic subgroup containing $X$; as a parabolic subgroup, $X_0$ is free abelian of rank at most $2$, and since it contains $X$ its rank must be exactly $2$. By assertion (3) of the proposition, which was proved above, $X_0$ is conjugate to the image of the inclusion homomorphism $\pi_1(T)\to\pi_1(M)$ for some torus component $T$ of $\partial K$. In particular, $X$ is conjugate to a subgroup of the image of the inclusion homomorphism $\pi_1(T)\to\pi_1(M)$. Since $K$ is a deformation retract of $M$ by Proposition \ref{lame song} and is therefore aspherical, it follows that $f$ is homotopic, in $K$, to a map of $T^2$ into $T$. This contradicts the definition of an essential singular torus.  \EndProof

\Number\label{octoberry}
An orientable hyperbolic $3$-manifold $M=\HH^3/\Gamma$ (or the Kleinian group $\Gamma$) is said to be {\it geometrically finite} if  $M$ is non-elementary and $\calc_1(M)$ has a compact truncation. 
\EndNumber

\Proposition\label{novemberry}
An orientable hyperbolic $3$-manifold $M$ has finite volume if and only if $M$ is geometrically finite and $\calc_1(M)$ has a truncation whose boundary components are all tori.
\EndProposition

\Proof
If  $\calc_1(M)$ has a truncation whose boundary components are all tori, it follows from Assertion (2) of Proposition \ref{lame song} that $\calc_1(M)$ has no boundary; hence $\calc_1(M)=M$. If in addition $M$ is geometrically finite, then $\calc_1(M)$ has finite volume according to \cite[Proposition 3.7]{M-T}; hence in this case $M$ has finite volume. Conversely, if $M$ has finite volume, then in particular $\calc_1(M)$ has finite volume, and \cite[Proposition 3.7]{M-T} implies that $M$ is geometrically finite. Furthermore, if we write $M=\HH^3/\Gamma$, the finiteness of $\vol M$ implies that the limit set of $\Gamma$ is the entire sphere at infinity; hence by \cite[Proposition 3.1]{M-T},  $\calc_1(M)$ has no boundary, and so the boundary components of a truncation of $\calc_1(M)$ are all tori.
\EndProof

An orientable hyperbolic $3$-manifold $M=\HH^3/\Gamma$ will be said to {\it have no rank-$1$ cusps} if every maximal parabolic subgroup of $\Gamma$ has rank $2$.

The following version of Thurston's geometrization theorem is a kind of converse to Proposition \ref{lame song}.

\Proposition\label{tho thad}
Let $K$ be a compact, irreducible orientable $3$-manifold which has non-empty boundary and has no essential singular tori. Then either $\pi_1(K)$ is isomorphic to either a Klein bottle group $\langle x,y:yxy^{-1}=x^{-1}\rangle$ or a free abelian of rank at most $2$, or $K$ is diffeomorphic to a truncation of $\calc_1(M)$ for some geometrically finite orientable hyperbolic $3$-manifold $M$ having no rank-$1$ cusps. 
\EndProposition

\Proof
Let $\calt$ denote the union of all torus components of $K$. If $\calt$ has a component $T$ such that the inclusion homomorphism $\pi_1(T)\to\pi_1(K)$ is not injective, then $\pi_1(K)$ is infinite cyclic by \cite[Proposition 3.10]{bounds}. We may therefore assume that $\pi_1(T)\to\pi_1(K)$ is injective for every component $T$ of $\calt$. We may also assume that $K$ is not a $3$-ball, as otherwise $\pi_1(K)$ is trivial.

Consider the case in which the pair $(K,\calt)$ has no essential singular annuli (\ref{erhard}). In this case, since $K$ also has no essential singular tori, and since $\pi_1(T)\to\pi_1(K)$ is injective for every component $T$ of $\calt$, the pair $(K,\calt)$ is a pared manifold in the sense of \cite[Definition 4.8]{morgan}. Furthermore, since $K$ is irreducible and is not a ball, and $\partial K\ne\emptyset$, the manifold $K$ is a Haken manifold in the sense defined on page 57 of \cite{morgan}. Hence the pared manifold $(K,\calt)$ satisfies the hypotheses of \cite[p. 70, Theorem B$'$]{morgan}. Since the components of $\calt$ (if any) are all tori (rather than annuli), the conclusion of \cite[p. 70, Theorem B$'$]{morgan} may be expressed, in the language of the present paper, by saying that there is a geometrically finite orientable hyperbolic manifold $M$, having no rank-$1$ cusps, such that $K$ is diffeomorphic to a truncation of $\calc_1(M)$.

There remains the case in which the pair $(K,\calt)$ has an essential singular annulus $f:(S^1\times[0,1],\partial(S^1\times[0,1]))\to (K,\calt)$. I will use the terminology of \cite{js} in the following argument. Let $(\Sigma,\Phi)$ denote a characteristic pair of the compact, irreducible pair $(K,\calt)$, which exists by \cite[p. 138]{js}. According to \cite[Remark IV.1.2]{js} and the definition of a characteristic pair (\cite[p. 138]{js}), after modifying $f$ by a homotopy of maps of pairs, we may assume that $f(S^1\times[0,1])\subset\Sigma$ and $f(\partial( S^1\times[0,1]))\subset\Phi$.  Let $\Sigma_0$ denote the component of $\Sigma$ containing $f(S^1\times[0,1])$.

According to the definition of a characteristic pair, $\Sigma_0$ is
a ``perfectly embedded pair'' in the sense of (\cite[p. 4]{js}).  In particular each component of $\Phi_0:=\Sigma_0\cap\partial M=\Sigma_0\cap\calt$ is a compact submanifold of $\calt$ whose boundary curves (if any) are homotopically non-trivial in $\calt$. Since the components of $\calt$ are tori, each component of $\Phi_0$ is a torus or an annulus.

According to the definition of a characteristic pair,  either $\Sigma_0$ is a Seifert fibered space and $\Phi_0\subset\Sigma_0$ is saturated, or $\Sigma_0$ is a $[0,1]$-bundle whose associated $\{0,1\}$-bundle is $\Phi_0$. In the latter subcase, since the components of $\Phi_0$ are tori and annulus, the base of the $[0,1]$-bundle is a torus, Klein bottle, annulus or M\"obius band, and hence the $[0,1]$-bundle may be given the structure of a Seifert fibered space in such a way that $\Phi_0$ is saturated. Thus in any event $\Sigma_0$ is a Seifert fibered space and $\Phi_0$ is saturated.

I claim that the components of $\Phi_0$ are tori. To prove this, let $\Sigma_1$ denote a regular neighborhood of the union of $\Sigma_0$ with all those components of $\calt$ which meet $\Sigma_0$. Since $\Phi_0$ is saturated in $\Sigma_0$, the manifold $\Sigma_1$ is a Seifert fibered space. Since $\Sigma_1\cap\calt$ is a union of boundary tori of $\Sigma_1$, the pair $(\Sigma_1,\Sigma_1\cap\calt)$ is a Seifert pair. Since $f$ is an essential singular annulus and $f(S^1\times[0,1])\subset\Sigma_0$, the inclusion map
$i:(\Sigma_1,\Sigma_1\cap\calt)\to(K,\calt)$ is an essential, nondegenerate map of Seifert pairs. The defining property of the characteristic pair therefore implies that $i$ is homotopic as a masp of pairs to a map $i'$ such that $i'(\Sigma_1)\subset\Sigma_0$ and $i'(\Sigma_1\cap\calt)\subset\Phi_0$. In particular, the inclusion of $\Sigma_1\cap\calt$ into $\calt$ is homotopic in $\calt$ to a map whose image is contained in $\Phi_0$. Since the components of $\Sigma_1\cap\calt$ are closed surfaces, it follows that $\Sigma_1\cap\calt=\Phi_0$, i.e. that the components of $\Phi_0$ are tori.

Now I claim that $K=\Sigma_0$. If this is false, we may fix a component $C$ of the frontier of $\Sigma_0$ in $K$. Since $\Sigma_0$ is a Seifert fibered space and $\Sigma\cap\partial K=\Phi_0$ is a union of components of $\partial\Sigma_0$, the surface $C$ is a torus. Since $\Sigma_0$ is perfectly embedded,
the inclusion homomorphism $\pi_1(C)\to\pi_1(M)$ is injective. As $K$ contains no essential tori, the inclusion $C\to K$ must be homotopic to a map of $C$ into $\partial K$. It then follows from \cite[Corollary 5.5]{waldhausen} that $C$ is boundary-parallel in $K$. This contradicts the definition of a perfectly embedded pair.

Since $K=\Sigma_0$, the manifold $K$ is a Seifert fibered space. Hence $\Gamma:=\pi_1(K)$ has a cyclic normal subgroup $N$ such that $Q:=\Gamma /N$ is a Fuchsian group. Let $\tGamma$ denote the centralizer of $N$ in $\Gamma$, so that $[\Gamma:\tGamma]\le2$. Let $q:\Gamma \to Q$ denote the quotient projection, and set $\tQ=q(\tGamma)$. Since $\partial K\ne\emptyset$, we may write $Q$ as a free product of non-trivial cyclic groups. 

If in the free product description of $Q$ there are at least three factors, or if there are at least two factors and one of them order strictly greater than $2$, then $\tQ$ has infinitely many non-conjugate maximal infinite cyclic subgroups. For each maximal cyclic subgroup $Z$ of $\tQ$, the group $q^{-1}(Q)$ is a maximal rank-$2$ free abelian subgroup of $\Gamma$. Hence $\Gamma$ has infinitely many non-conjugate maximal rank-$2$ free abelian subgroups. As $\partial K$ has only finitely many components, it follows that $K$ admits an essential singular torus, a contradiction to the hypothesis.

Hence in the free product description of $Q$ there are at most two factors, and if there are two factors they are both of order $2$. This means that either the base $B$ of the Seifert fibration is a M\"obius band or annulus and there are no singular fibers; or $B$ is a disk and there is at most one singular fiber; or $B$ is a disk and there are two singular fibers, both of local degree $2$. It follows that $K$ is diffeomorphic to a solid torus---a contradiction---or to $T^2\times[0,1]$ or a twisted $I$-bundle over a Klein bottle. Hence $\pi_1(K)$ is either a Klein bottle group or a free abelian group of rank $2$.
\EndProof

\section{Geometrically finite Kleinian groups}\label{GF section}

The main results of the section are Proposition \ref{no great loss}, which I have not seen stated before, and Proposition \ref{onion omelette}, which is routine. The following elementary result, Proposition \ref{i thought i was dead}, which is needed for the proof of  \ref{no great loss}, seems surprisingly difficult to locate in the literature.

\Proposition\label{i thought i was dead}
Let $V\subset\CC^N$ be an affine algebraic set defined over an algebraically closed subfield $F$ of $\CC$. Then $V\cap F^N$ is a dense subset of $V$ in the classical (complex) topology.
\EndProposition

\Proof
According to \cite[Theorem 3.10.9]{nagata}, if an algebraic set is defined over the algebraically closed subfield $F$ of $\CC$, then its irreducible components are defined over $F$.

I will prove the assertion by induction on the dimension of $V$, which is defined to be the maximum of the dimensions of its irreducible components. If $\dim V=0$ then $V=\{P_1,\cdots,P_k\}$ is a finite set. Since each irreducible component $\{P_i\}$  of $V$  is defined over $F$, it follows for example from \cite[Corollary 30.3]{isaacs} that $P_i\in F^{N}$, which gives the conclusion in this case.

Now suppose that $V\subset\CC^{N}$ is an affine algebraic set defined over $F$ and having dimension $d>0$, and that the assertion is true for affine algebraic subsets of $\CC^N$ defined over $F$ and having dimension less than $d$. We must show that $V\cap F^N$ is dense in the classical (complex) topology of  $V$. Since the irreducible components of $V$ are defined over $F$, we may assume without loss of generality that $V$ is irreducible.

Let $U$ be a non-empty subset of $V$ which is open in the classical topology. We must show that $U\cap{F}^{N}\ne\emptyset$. Let $P_0$ be a point of $U$. Since $\dim V>0$, there is an irreducible curve $X$ in $\CC^{N}$ with $P\in X\subset V$. Let $X_0\subset X$ denote the set of smooth points of $X$. Then $X_0$ is non-empty and Zariski-open in $X$, and by \cite[Theorem 2.33]{mumford} it is a dense subset of $X$ in the classical topology. Hence $X_0\cap U\ne\emptyset$. Fix a (classically) connected component $W$ of $X_0\cap U$. Then $W$ is a complex $1$-manifold; in particular it contains more than one point, and hence one of the coordinate functions on $\CC^{N}$ is non-constant on $W$. If $c$ denotes such a coordinate function, then $c|W$ is a non-constant holomorphic function on the connected complex $1$-manifold $W$, and is therefore an open map to $\CC$. Since $F$ is dense in $\CC$, it follows that $c(W)\cap F\ne\emptyset$. Choose a point $\alpha\in c(W)\cap F$. If we set $V'=c^{-1}(\alpha)$, it follows that $U':=V'\cap U\ne\emptyset$. But $V'$ is an algebraic set defined over $F$. It is a proper subset of $V$ since $c|W$ is non-constant, and since $V$ is irreducible it follows that $\dim V'<d$. By the induction hypothesis, $V'\cap F^N$ is dense in the classical topology of $V'$, and hence $U'\cap{F}^{N}\ne\emptyset$. In particular
$U\cap{F}^{N}\ne\emptyset$.
\EndProof

\Proposition\label{no great loss}
Let $\Gamma$ be a non-elementary, finitely generated, torsion-free Kleinian group. Then $\Gamma$ is isomorphic to a Kleinian group $\Gamma_1$ such that (i) $\Gamma_1$ is geometrically finite, (ii) every maximal parabolic subgroup of $\Gamma_1$ has rank $2$, and (iii) $\Gamma_1\le\pizzle(E)$ for some number field $E$.
\EndProposition

\Proof
If $M$ denotes the orientable hyperbolic $3$-manifold $\HH^3/\Gamma$, we have $\Gamma\cong\pi_1(M)$. I claim:
\Claim\label{make mine finite}
$\Gamma$ is isomorphic to $\pi_1(M_0)$ for some geometrically finite orientable hyperbolic $3$-manifold $M_0$ having no rank-$1$ cusps.
\EndClaim
Indeed, \ref{make mine finite} is obvious in the case where $M$ is closed, since we may then take $M_0=M$. If $M$ is not closed, let us fix a truncation $K$ of $\calc_1(M)$, and denote by $\calt$ the union of all torus components of $\partial K$. According to the main theorem of \cite{sullivan}, $\Gamma$ has only finitely many conjugacy class of maximal parabolic subgroups; hence by Proposition \ref{lame song}, $\calt$ has only finitely many components, i.e. it is compact. It then follows from the Relative Core Theorem of \cite{mcc} that there is a compact, connected submanifold $K_0\supset \calt$ of $K$ such that the inclusion homomorphism $\pi_1(K_0)\to\pi_1(K)$ is an isomorphism. According to \cite[Lemma 3.4]{bounds}, there is a compact, irreducible, $3$-dimensional submanifold $K_1$ of $M$ such that $K_1\supset K_0$, and such that the inclusion homomorphism $\pi_1(K_0)\to\pi_1(K_1)$ is surjective. It follows that
$K_1\supset \calt$ and that the inclusion homomorphism $\pi_1(K_1)\to\pi_1(K)$ is an isomorphism. Since $\pi_1(K_1)\cong\Gamma$ is torsion-free and non-abelian, and hence infinite, it follows from \cite[Proposition 3.9]{bounds} that $K_1$ is aspherical; on the other hand, $K$ is homotopy equivalent to the hyperbolic manifold $M$ and is therefore aspherical. Hence the inclusion $K_1\to K$ is a homotopy equivalence, and so $K_1$ is a strong deformation retract of $K$. Since $K$ has no essential tori by Proposition \ref{lame song}, any singular essential torus $f:T^2\to K_1$ would be homotopic in $K$ to a map $g$ of $T^2$ into $\calt\subset K_1$; since $K_1$ is a strong deformation retract of $K$, the maps $f$ and $g$ would be homotopic in $K_1$, a contradiction to the essentiality of $f$. Hence $K_1$ has no essential tori. We have $\partial K_1\ne\emptyset$ since $M$ is non-compact. Since $\pi_1(K)$ is isomorphic to the non-elementary Kleinian group $\Gamma$, it has no abelian subgroup of finite index. It therefore follows from Proposition \ref{tho thad} that $K_1$ is diffeomorphic to $\calc_1(M_0)$ for some geometrically finite, orientable hyperbolic $3$-manifold $M_0$ having no rank-$1$ cusps. In particular we have $\Gamma\cong\pi_1(M)\cong\pi_1(K_1)\cong\pi_1(M_0)$. Thus \ref{make mine finite} is established in all cases.

Let $\Hom(\Gamma,\zzle(\CC))$ denote the set of all representations of $\Gamma$ in $\zzle(\CC)$. Let $\Hom(\Gamma,\zzle(\CC))^*$ denote the subset of $\Hom(\Gamma,\zzle(\CC))$ consisting of all representations $\rho$ such that $\rho$ maps each parabolic element of $\Gamma$ to an element of trace $\pm2$ in $\zzle(\CC)$. Fix a generating set $x_1,\ldots,x_t$ of $\Gamma$, and define a map of sets $\Phi:\Hom(\Gamma,\zzle(\CC))\to\zzle(\CC)^t\subset(\calM_2)^t$ by $\Phi(\rho)=(\rho(x_i))_{1\le i\le t}$. Then $R:=\Phi(\Hom(\Gamma,\zzle(\CC))\subset(\calM_2)^t$  is readily seen to be an affine algebraic subset of $(\calM_2)^t$ (see \cite[Subsection 4.1]{handbook}). Furthermore, $R^*:=\Phi( \Hom(\Gamma,\zzle(\CC)^*)\subset R$ is also an algebraic set. (In the notation of \cite[Subsection 4.4]{handbook},  $R^*$ is the locus of zeros within $R$ of the polynomials $I_\gamma^2-4$, where $\gamma$ ranges over the parabolics in $\Gamma$.

Now let $P:\zzle(\CC)\to\pizzle(\CC)$ denote the quotient map. Let $\calf\calg$ denote the set of all points of $R$ of the form $\Phi(\rho)$, where $\rho:\Gamma\to\zzle(\CC)$  is a representation of $\Gamma$ such that $P\circ\rho$ is faithful and $P\circ\rho(\Gamma)$ is a geometrically finite Kleinian group of which all maximal parabolic subgroups have rank $2$. If $M_0$ is the hyperbolic manifold given by \ref{make mine finite}, then according to \cite[Proposition 3.1.1]{splittings}, the discrete faithful representation of $\pi_1(M_0)$ in $\pizzle(\CC)$ defined by the hyperbolic structure of $M_0$ may be lifted to a representation $r:\pi_1(M_0)\to\zzle(\CC)$; precomposing $r$ with an isomorphism of $\Gamma$ onto $\pi_1(M_0)$ gives a representation $\rho$ such that $\Phi(\rho)\in\calf\calg$. Hence:
\Claim\label{not for nothing} $\calf\calg\ne\emptyset$. 
\EndClaim

Since a discrete  rank-$2$ free abelian subgroup of $\zzle(\CC)$ must be parabolic, we have $\calf\calg\subset R^*$. Now I claim:
\Claim\label{situation normal all fogged up}
The subset $\calf\calg$ is open in $R^*$.
\EndClaim

To prove \ref{situation normal all fogged up} one must show, given a point $\Phi(\rho_0)\in\calf\calg$, that some neighborhood of $\Phi(\rho_0)$ in $R^*$ is contained in $\calf\calg$. Here $P\circ\rho_0$ is faithful, $\Gamma_0:=P(\rho_0(\Gamma))$ is discrete and geometrically finite, and  all its  maximal parabolic subgroups have rank $2$. Since the Kleinian group $\Gamma_0$ is geometrically finite, it follows from \cite[Theorem 3.7]{M-T} that there is a finite-sided Dirichlet polyhedron for $\Gamma_0$.  It then follows from the proof of \cite[Proposition 9.2]{marden} that $\Phi(\rho_0)$ has a neighborhood in $R^*$ consisting of points of the form $\Phi(\rho)$, where $\rho$ is a faithful representation such that $P\circ\rho(\Gamma)$ is discrete, and hence $\Phi(\rho)\in\calf\calg$. This completes the proof of \ref{situation normal all fogged up}.

(The definition of a Kleinian group used in  \cite{marden} includes the condition that the group have a non-empty set of discontinuity on the sphere at infinity, although this condition does not appear to be used in the proof of  \cite[Proposition 9.2]{marden}. This is why, in the proof of \ref{situation normal all fogged up} given above, I have had to quote the proof of  \cite[Proposition 9.2]{marden} rather than the statement; the latter result, {\it as stated}, does not cover the case in which the set of discontinuity $\Omega_0\subset S_\infty$ of $\Gamma_0$ is empty. An alternative approach in the case $\Omega_0=\emptyset$ is to apply \cite[Proposition 3.1]{M-T} to the manifold $M_0=\HH^3/\Gamma_0$ to deduce that $\calc_1(M_0)$ has empty boundary. It then follows from the definition of a truncation that $\calc_1(M_0)$ has a truncation whose boundary components are all tori. Hence by Proposition \ref{novemberry}, $M$ has finite volume. In the finite-volume case, \ref{situation normal all fogged up} is a well-known consequence of the results of \cite{garland}.)

Let $\caly$ denote the set of all points of $R^*$ whose coordinates are algebraic numbers. Since $\caly$ is the locus of zeros of a set of polynomial equations with integer coefficients, we may apply Proposition \ref{i thought i was dead}, taking $F$ to be the algebraic closure of $\QQ$ in $\CC$, to deduce that $\caly$ is dense in $R^*$. In view of \ref{not for nothing} and \ref{situation normal all fogged up} it follows that $\caly\cap\calf\calg\ne\emptyset$. Let $\rho_1$ be a representation such that $\Phi(\rho_1)\in\caly\cap\calf\calg$, and set $\Gamma_1=P(\rho_1(\Gamma))$ under the quotient homomorphism $\zzle(\CC)\to\pizzle(\CC)$. Since 
$\Phi(\rho_1)\in\calf\calg$, the group $\Gamma_1$ is discrete, is isomorphic to $\Gamma$, and satisifies conditions (i) and (ii) of the conclusion of the proposition.
Since $\Phi(\rho_1)\in\caly$, there is a number field $E$ containing the coordinates of $\Phi(\rho_1)$. It follows that $\rho_1(\Gamma)\le\zzle(E)$, and condition 
(iii) of the conclusion follows.
\EndProof

\Proposition\label{onion omelette}
Let $\Gamma$ be a  geometrically finite, torsion-free Kleinian group, and let $R$ be a positive real number. Then the set of loxodromic elements of $\Gamma$ having length less than $R$ is a union of finitely many conjugacy classes.
\EndProposition

\Proof
Set $M:=\HH^3/\Gamma$. Fix a truncation $K$ of $\calc_1(M)$, and let $X_1,\ldots,X_k$ denote the standard cusp neighborhoods bounded by the torus components of $\partial K$. Since $\Gamma$ is geometrically finite, $K$ is compact (cf. \ref{vexcor}). Let $D$ denote the intrinisic diameter of $K$, so that any two points of $K$ are joined by a path of length at most $D$. Fix a base point $\star\in K$, and fix a point $\widetilde\star\in\HH^3$ that maps to $\star$ under the quotient map $\HH^3\to M$. The point $\widetilde\star$ determines an isomorphism $J:\Gamma\to\pi_1(M,\star)$.

Let $\gamma \in\Gamma$ be a  loxodromic element of length at most $R$. Then the conjugacy class of $J(\gamma )\in\pi_1(M,\star)$ is represented by an oriented closed geodesic $C$ in $M$ with length at most $R$. In the notation of \ref{vexcor} we have $|C|\subset\calc(M) \subset\calc_1(M)$. Since $\gamma $ is loxodromic, $|C|$ cannot be contained in any of the $X_i$. Hence $|C|\cap K\ne\emptyset$. Let us fix a point $y\in S^1$ such that $C(y)\in K$. Let $\beta$ be a positively oriented loop in $S^1$, based at $y$, which defines a generator of $\pi_1(S^1,y)$; then $c:=C\circ\beta$ is a loop in $M$ based at $\beta(y)$.  Let $\alpha$ be a path in $K$ which begins at $\star$, ends at $\beta(y)$, and has length at most $D$. Then $c'=\alpha\star c\star\bar\alpha$ is
a loop based at $\star$ which has length at most $R+2D$, and $[c']$ is conjugate to $J(\gamma )$ in $\pi_1(M,\star)$. Hence $\gamma '=J^{-1}([c'])$ is conjugate to $\gamma $ in $\Gamma$, and $d(\widetilde\star,\gamma \cdot\widetilde\star)\le R+2D$. Since $\Gamma$ is discrete there are only finitely many elements of $\Gamma$ that displace $\widetilde\star$ by a distance at most $R+2D$, and the conclusion follows.
\EndProof

\Corollary\label{cat}
Let $\tGamma$ be a subgroup of $\zzle(\CC)$ which maps isomorphically onto a geometrically finite, torsion-free Kleinian group under the quotient homomorphism $\zzle(\CC)\to\pizzle(\CC)$. Let $\tau\ne\pm2$ be a  complex number. Then the set of elements of $\tGamma$ having trace $\tau$ is a union of finitely many conjugacy classes. 
\EndCorollary

\Proof
Let $\Gamma$ denote the image of $\tGamma$ on $\pizzle(\CC)$.
If $\gamma\in\tGamma$ has trace $\tau$, then $[\gamma]\in\Gamma$ is loxodromic since $\tau\ne\pm2$; and if $l\in(0,\infty)$ and $\theta\in\RR/2\pi\RR$ denote the length and twist angle of $\gamma$, then the quantity $2\cosh((l+i\theta)/2)$, which is well-defined up to sign, is equal to $\pm\tau$. Hence all elements of $\tGamma$ with trace $\tau$ map to loxodromic elements of the same length in $\Gamma$, and the assertion of the corollary follows from Proposition \ref{onion omelette}.
\EndProof

\section{Dehn filling}\label{filling section}

The main result of this section is Proposition \ref{de uddah joik}, which I discussed in the introduction; as I pointed out there, it builds on ideas from \cite{AllHam}. Lemma \ref{some kina wonnaful} will be needed for the proof.

\Lemma\label{some kina wonnaful}
Let $K$ be a compact irreducible orientable $3$-manifold with no essential singular tori. Let $T$ be a torus component of $\partial K$, and let $x$ denote a base point in $T$. Suppose that the inclusion homomorphism  $\pi_1(T,x)\to\pi_1(M,x)$ is injective, and let $P$ denote its image. Then there
 exist a compact, irreducible orientable $3$-manifold $K'$, a component $T'$ of $\partial K'$, a base point $x'\in T'$ and a homomorphism $J:\pi_1(K,x)\to \pi_1(K',x')$, such that 
\begin{enumerate}
\item every component of $\partial K'$ is a torus;
\item there are no essential tori in $K'$; and
\item $J|P$ is an isomorphism of $P$ onto the image of the inclusion homomorphism  $\pi_1(T',x')\to\pi_1(K',x')$.
\end{enumerate}
\EndLemma

\Proof 
Let
$T_1,\ldots,T_k$ denote the torus components of $\partial K$, indexed so that $T_1=T$.
If $T_1,\ldots,T_k$ are the only components of $\partial K$, the conclusion of the lemma follows upon setting $K'=K$ and taking $J$ to be the identity map. I will therefore assume that  $\partial K$ has one or more boundary component of genus greater than $1$; let $F_1,\ldots,F_n$ denote the higher-genus components of $\partial \calc_1(M)$, and let $g_j>1$ denote the genus of $F_i$.

First consider the case in which $K$ is boundary-irreducible. For $j=1,\ldots,n$, the construction of \cite{fujii} gives a compact orientable hyperbolic $3$-manifold $Q_j$ with connected totally geodesic boundary such that $\partial Q_j$ has genus $g_j$. Since $Q_i$ has totally geodesic boundary, it follows from \cite[Theorem 2.1]{bonahon} that $Q_i$ is irreducible and boundary-irreducible, $Q_i$ has no essential singular tori and. 
$(Q_i,\partial Q_i)$ has no essential singular 
annuli (see \ref{erhard}).

Let $K'$ denote the orientable $3$-manifold obtained from the disjoint union of $K$ and $Q_1,\ldots,Q_n$ by gluing $F_j$ to $\partial Q_j$ via some (arbitrarily chosen) diffeomorphism $\eta_j: F_j\to\partial Q_j$ for $j=1,\ldots,n$. Since $K$ and each $Q_i$ are boundary-irreducible, the surface $F=F_1\cup\dots\cup F_n$ is incompressible in $K'$. 
I claim that $K'$ is irreducible and boundary-irreducible, and has no essential singular tori.

To prove irreducibility, suppose that $\Sigma\subset K'$ is a $2$-sphere. After a small isotopy we may assume that $\Sigma$ meets $F$ transversally. I will prove by induction on the number $c$ of components of $\Sigma\cap F$ that $\Sigma$ bounds a ball in $K'$. The case $c=0$ follows from the irreducibility of $K$ and the $Q_j$. If $c>0$, each component of  $\Sigma\cap F$ bounds a disk in  $F$ since $F$ is incompressible. Among all disks in $F$ bounded by components of  $\Sigma\cap F$, choose one, say $D$, which is minimal with respect to inclusion. Then $C:=\partial D$ bounds two disks $E_1,E_2\subset\Sigma$, and each $E_i\cup D$ is a $2$-sphere $\Sigma'_i$ which is isotopic by a small isotopy to a $2$-sphere meeting $F$ in fewer than $c$ components. By the induction hypothesis, each $\Sigma_i'$ bounds a ball $B_i$. We have either $B_1\cap B_2=D$ or, after possibly re-indexing the $B_i$, that $B_1\subset B_2$. In the first case $B_1\cap B_2$ is a ball bounded by $\Sigma$, and in the second case $\overline{B_2- B_1}$ is a ball bounded by $\Sigma$.

The proof of boundary-irreducibility is somewhat similar. Suppose that $\Delta\subset K'$ is a properly embedded disk which meets $F$ transversally. I will prove by induction on the number $c$ of components of $\Delta\cap F$ that $\partial\Delta$ bounds a disk in $\partial K'$. The case $c=0$ follows from the boundary-irreducibility of $K$ and the $Q_j$. If $c>0$, each component of  $\Delta\cap F$ bounds a disk in  $F$ since $F$ is incompressible. Among all disks in $F$ bounded by components of  $\Delta\cap F$, choose one, say $D$, which is minimal with respect to inclusion. Then $C:=\partial D$ bounds a disk $E\subset\Delta$, and $(\Delta-E)\cup D$ is a disk $\Delta'$ which is isotopic by a small isotopy to a disk meeting $F$ in fewer than $c$ components. By the induction hypothesis, $\partial\Delta=\partial\Delta'$ bounds a disk in $\partial K'$.

To show that $K'$ has no essential singular tori, suppose that $f:T^2\to K$ induces an injection of fundamental groups. After a small homotopy we may assume that $f$ is transverse to $F$. I will prove by induction on the number $c$ of components of $f^{-1}(F)$ that $f$ is homotopic to a map into $\partial K'$. The case $c=0$ follows from the fact that $K$ and the $Q_j$ have no essential singular tori. Now suppose that $c>0$ and that some component $\gamma$ of $f^{-1}(F)$ bounds a disk $\Delta\subset T^2$. Since $F$ is incompressible, $f|\gamma$ is homotopically trivial in $F$. Hence there is a map $f':T^2\to K'$ which agrees with $f$ outside $\Delta$ and maps $\Delta$ into $F$. Since $K'$ is irreducible we have $\pi_2(K')=0$, and hence $f$ is homotopic to $f'$. Clearly
$f'$ is in turn homotopic by a small homotopy to a map $f_1'$ such that $(f_1')^{-1}(F)$ has fewer than $c$ components. By the induction hypothesis $f_1'$ is homotopic to a map into $\partial K'$, and hence so is $f$.

Now suppose that $c>0$ and that no component $\gamma$ of $f^{-1}(F)$ bounds a disk $\Delta\subset T^2$. Since $c>0$ we have $f^{-1}(Q_1\cup\cdots Q_n)\ne\emptyset$. After re-indexing the $Q_i$ we may assume that $f^{-1}(Q_1)\ne\emptyset$. Choose a component $A$ of $f^{-1}(Q_1)$. Since $(Q_1,\partial Q_1)$ has no essential annuli, $f|A$ is homotopic in $Q_1$ to a map of $A$ into $F_1$. 
Hence $f$ is homotopic to a map $f':T^2\to K'$ which agrees with $f$ outside $A$ and maps $A$ into $F$. Clearly
$f'$ is in turn homotopic by a small homotopy to a map $f_1'$ such that $(f_1')^{-1}(F)$ has fewer than $c$ components. By the induction hypothesis $f_1'$ is homotopic to a map into $\partial K'$, and hence so is $f$.

Now let $J:\pi_1(K,x)\to \pi_1(K',x')$ denote the inclusion homomorphism. I claim the conclusions of the lemma hold with this choice of $J$ if we set $T=T'$ and $x'=x$. By construction $K'$ is compact and orientable and its boundary components are tori. I have shown that $K$ is irreducible and contains no essential singular tori. It remains only to observe that Conclusion (iii) of the lemma holds. If $P'$ denotes the image of the inclusion homomorphism  $\pi_1(T,x)\to\pi_1(K',x')$, the definitions of $J$, $P$ and $P'$ imply that $J(P)=P'$. That $J|P:P\to P'$ is an isomorphism is tantamount to saying that the inclusion homomorphism $\pi_1(T,x)\to\pi_1(K',x')$ is injective; this follows from the boundary-irreducibility of $K'$, which I proved above.

This completes the proof of the lemma in the case where $K$ is boundary-irreducible.

Now consider the case in which $\partial K$ is  boundary-reducible. Let $D_1,\ldots,D_m$ be a maximal system of pairwise disjoint, non-boundary-parallel, pairwise non-parallel, properly embedded disks in $K$. Let $E$ be a regular neighborhood of $D_1\cup\cdots\cup D_m$ in $K$, and set $K_0=\overline{K-E}$. 

Since
the inclusion homomorphism $\pi_1(T)\to\pi_1(K)$ is injective, none of the disks $D_1,\ldots,D_m$ has its boundary in $T$, and hence $T\subset K_0$.
Let $K_1$ denote the component of $K_0$ that contains $T$.

Since $K$ is irreducible and has no essential singular tori, and since $K_1$ is a component of a manifold obtained from $K$ by splitting it along a collection of pairwise disjoint properly embedded disks, $K_1$ is itself irreducible and has no essential singular tori. Furthermore, $T$ is a component of $\partial K_1$ and the inclusion $\pi_1(T,x)\to\pi_1(K_1,x)$ is injective; I will denote its image by $P_1$. On the other hand, the maximality of the family $D_1,\ldots,D_m$ implies that $K_1$ is boundary-irreducible. It therefore follows from the case of the lemma already proved that there 
 exist a compact, irreducible, orientable $3$-manifold
$K'$, a component $T'$ of $\partial K'$, a base point $x'\in T'$ and a homomorphism $J_1:\pi_1(K_1,x)\to \pi_1(K',x')$ such that conclusions (i)---(iii) of the lemma hold with $K_1$, $P_1$ and $J_1$  in place of $K$, $P$ and $J$.

Since $K_1$ is a component of a manifold obtained from $K$ by splitting along a family of disks, the inclusion homomorphism $I:\pi_1(K_1,x)\to\pi_1(K,x)$ maps $\pi_1(K_1,x)$ isomorphically onto a free factor of $\pi_1(K,x)$. In particular $I$ has a left inverse $r$ which is a homomorphism. The definitions of $I$, $P_1$ and $P$ imply that $I(P_1)=P$; since $I$ is injective, it restricts to an isomorphism of  $P_1$ onto $P$. Hence $r$ maps $P$ isomorphically onto $P_1$. Since $J_1$ maps $P_1$ isomorphically onto $P'$, the homomorphism $J:=J_1\circ r:\pi_1(K,x)\to \pi_1(K',x')$ maps $P$ isomorphically onto $P'$, and the lemma is proved in this  case as well.
\EndProof

\Proposition\label{de uddah joik}
Let $\Gamma$ be a torsion-free geometrically finite  Kleinian group such that every maximal parabolic subgroup of $\Gamma$ has rank $2$. Let $P$ be a maximal parabolic subgroup of $\Gamma$. Then there is a finite set $Y\subset P$ such that for every $\gamma\in P-Y$, and every positive integer $m$, there
 exist a Kleinian group $\hatGamma$ and a homomorphism $H:\Gamma\to \hatGamma$ such that $H(\gamma)$ has order $m$ in $\hatGamma$.
Furthermore, if $Y$ is such a set, then every positive integer is a finitistic order for every $\gamma\in P-Y$.
\EndProposition

\Proof

Set $M=\HH^3/\Gamma$, and fix a truncation $K$ of $\calc_1(M)$. Let $T_1\cup\cdots\cup T_k$ denote the torus boundary components of $K$. By Proposition \ref{lame song}, $K$ is. Since $P$ is a maximal parabolic subgroup of $\Gamma$, it follows from Proposition \ref{lame song} that, after possibly re-indexing  the $T_i$ and choosing a base point $x\in T_1$, we have an isomorphic identification of $\Gamma$ with $\pi_1(K,x)$ under which $P$ is the image of the inclusion homomorphism $\pi_1(T_1,x)\to \pi_1(K,x)$.

I will begin by proving the first assertion in the case where $\Gamma$ has finite covolume. In this case we have $\calc_1(M)=M$; and if $K$ is a truncation of $M$, then $\partial K$ is a union of tori $T_1,\ldots, T_k$ by Proposition \ref{novemberry}. The proof in this case is an application of  Thurston's hyperbolic Dehn filling theorem, which is stated and proved as \cite[Theorem 2.1]{pe-po}.

The following notation will be borrowed from \cite{pe-po}. Fix a basis $\lambda_i,\mu_i$ of $H_1(T_i;\ZZ)$ for each $i$. Denote by $C$ the set
of all coprime pairs of integers, together with a symbol $\infty$. For
$c_1,\dots,c_k\in C$ denote by $M_{c_1\cdots c_k}$ the manifold obtained from
$K$ as follows: if $c_i=\infty$, glue a half-open collar $T_i\times[0,\infty)$ to $T_i$; if $c_i=(p_i,q_i)$, glue a solid torus $J_i=D^2\times S^1$ to
$K$ along $T_i$, by gluing the meridian
$S^1\times\{1\}$ to a curve representing the class $p_i\lambda_i+q_i\mu_i$. Consider the set
$$G=\{\infty\}\cup\{g\in\RR^2:\ g=r\cdot(p,q){\rm\ for\ some\ }
r>0{\rm\ and\ relatively\ prime\ }p,q\in\ZZ\}.$$
For $g=r\cdot(p,q)\in G\setminus\{\infty\}$
define $c(g)=(p,q)$ and
$\vartheta(g)=2\pi/r$. Set $c(\infty)=\infty$. Topologize $G$ as a subset of
$\RR^2\cup\{\infty\}=S^2$. 

According to \cite[Theorem 2.1]{pe-po}
there exists a neighbourhood
$\calf$ of $(\infty,\dots,\infty)$ in $G^k$ such that for
$(g_1,\dots,g_k)\in\calf$ the manifold $M_{c(g_1)\dots c(g_k)}$ admits the
structure of a complete finite-volume hyperbolic cone manifold in which the cone
locus consists of the cores $\{0\}\times S^1$ of those 
$J_i$ such that $g_i\neq\infty$, and the cone angle at the core of $J_i$ is $\vartheta(g_i)$.

The inclusion homomorphism $\pi_1(T_1,x)\to\pi_1(K,x)$ may be regarded as an isomorphism of $\pi_1(T_1,x)$ onto $P$. Composing the inverse of this isomorphism with the Hurewicz isomorphism from $\pi_1(T_1,x)$ to $H_1(T_1,\ZZ)$, one obtains an isomorphism $e:P\to H_1(T_1,\ZZ)$. The free abelian group $P$ has a basis consisting of the elements ${\frak l}:=e^{-1}(\lambda_1)$  and ${\frak m}:=e^{-1}(\mu_1)$.

Now choose an integer $B>0$ such that for every $c=(x,y)\in\RR^2$ with $\max(|x|,|y|)>B$ we have $(c,\infty,\ldots,\infty)\in\calf$. 
Let $Y$ denote the finite subset of $P$ consisting of all elements of the form ${\frak l}^r{\frak m}^s$ as $r$ and $s$ range over all integers of absolute value at most $B$.
I will show that the first assertion holds with this choice of $Y$.

Let $\gamma\in P-Y$ be given, and let  $m$ be a positive integer. Let us write $\gamma={\frak l}^{dp}{\frak m}^{dq}$, where $p$ and  $q$ are relatively prime integers,
$d$ is a positive integer,  and $\max(|dp|,|dq|)>B$. Hence $\max(md|p|,md|q|)>B$. If we set $g=(mdp,mdq)$, it follows that $(g,\infty,\ldots,\infty)\in\calf$. By definition we have $c(g)=(p,q)$ and $\vartheta(g)=2\pi/(md)$.

The manifold $M_{c(g),\infty,\ldots,\infty}$ is obtained from $ K$ by attaching a solid torus $J$ along the boundary component $T_1$, and attaching half-open collars to the remaining components of $\partial K$. Here $J$ may be given a product structure $J=D^2\times S^1$ in such a way that $S^1\times\{1\}$ is glued to a simple closed curve representing the class $p\lambda+q\mu\in H_1(K;\ZZ)$. 
The defining property of $\calf$ implies that $M_{c(g),\infty,\ldots,\infty}$ admits the structure of a complete finite-volume hyperbolic cone manifold, whose cone locus consists of the core $\{0\}\times S^1$ of $J$, and the cone angle at this core is $\vartheta(g) =2\pi/(md)$. Since $md$ is an integer, the cone manifold structure gives $M_{c(g),\infty,\ldots,\infty}$ the structure of a hyperbolic orbifold whose singular locus is the core curve; in the neighborhood of any point of the singular locus, the orbifold is topologically modeled on the quotient of $\RR^3$ by the finite cyclic group $\ZZ/md$, acting through rotations about an axis. It follows that the orbifold fundamental
group $\pi_1^{\rm orb}(M_{c(g),\infty,\ldots,\infty})$ is isomorphic to a Kleinian group, and that the inclusion homomorphism (defined modulo inner automorphisms) from $\pi_1(K)$ to $\pi_1^{\rm orb}(M_{c(g),\infty,\ldots,\infty})$ maps $\gamma_0={\frak l}^{p}{\frak m}^q$ onto a element of order of $md$. Hence $H$ maps $\gamma=\gamma_0^d$ onto an element of order $m$. This proves the first assertion of the proposition in the finite-covolume case.

I will now prove the first assertion in the general case. According to Proposition \ref{lame song}, $K$ is irreducible and has no essential singular tori. The inclusion homomorphism  $\pi_1(T,x)\to\pi_1(M,x)$ is injective by Proposition \ref{lame song}, and it has image $P$ according to our choice of $T$. Hence we may apply Lemma \ref{some kina wonnaful}, taking $T=T_1$. This gives a compact, irreducible orientable $3$-manifold $K'$,  a component $T'$ of $\partial K'$, a base point $x'\in T'$ and a homomorphism $J:\pi_1(K,x)\to \pi_1(K',x')$. Furthermore, every component of $\partial K'$ is a torus; there are no essential tori in $K'$; and if $P'$ denotes the image of the inclusion homomorphism  $\pi_1(T,x')\to\pi_1(K',x')$, then $J$ restricts to an isomorphism of $P$ onto $P'$. 

I claim:

\Claim\label{topo giggio}
There
is a finite set $Y'\subset P'$ such that for every $\gamma'\in P'-Y'$, and every positive integer $m$,  there
 exist a Kleinian group $\hatGamma$ and a homomorphism $h':\pi_1(K',x')\to \hatGamma$ such that $h'(\gamma')$ has order $m$ in $\hatGamma$. 
\EndClaim

To prove (\ref{topo giggio}), first note that by Proposition \ref{tho thad}, either (i) $K'$ is diffeomorphic to a truncation of $\calc_1(M')$ for some geometrically finite orientable hyperbolic $3$-manifold $M'$ having no rank-$1$ cusps, (ii) $\pi_1(K',x')$ is free abelian, or (iii) $\pi_1(K',x')$ is a Klein bottle group. 

If (i) holds, then since the boundary components of $K'$ are all tori, $M'$ has finite volume by Proposition \ref{novemberry}, and we may write $M'=\HH^3/\Gamma'$ for some torsion-free Kleinian group $\Gamma'$ of finite covolume. Furthermore,  we may identify $\pi_1(K',x')$ isomorphically with $\Gamma'$ in such a way that $P'$ is a maximal parabolic subgroup. In this case, (\ref{topo giggio}) follows from the finite-covolume case of the first assertion of the lemma, which has already been proved.

If (ii) or (iii) holds, I will set $\Gamma'=\pi_1(K',x')$. If (ii) holds, then according to the construction of \ref{uneeda}, for every $\gamma'\in\Gamma'-\{1\}$ and every positive integer $m$, there is a homomorphism $h'$ of $\Gamma'$ onto a finite cyclic group $Z$ such that $h'(\gamma')$ has order $m$; as $Z$ is in particular isomorphic to a Kleinian group, (\ref{topo giggio}) holds in this case with $Y'=\{1\}$. 

If (iii) holds, we may identify $\Gamma'$ isomorphically with $\langle x,y: yxy^{-1}=x^{-1}\rangle$. Suppose that $\gamma'\in \Gamma'-\{1\}$ and a positive integer $m$ are given. If $\gamma'$ is not a power of $x$, its image $\gamma''$ under the quotient homomorphism $\Gamma'\to\Gamma'/\langle\langle x\rangle\rangle$ is non-trivial. Since $\Gamma'/\langle\langle x\rangle\rangle$ is infinite cyclic, the construction of \ref{uneeda} gives a homomorphism $h''$ of $\Gamma'/\langle\langle x\rangle\rangle$ onto a finite cyclic group $Z$ such that $h''(\gamma'')$ has order $m$, and again (\ref{topo giggio}) holds with $Y'=\{1\}$. If $\gamma'=x^k$ for some $k\ne0$, and if $D$ denotes the finite dihedral group $\langle u,v:v^2=1,u^{m|k|}=1,vuv^{-1}=u^{-1}\rangle$, the homomorphism $h':\Gamma'\to D$ defined by $h'(x)=u$, $h'(y)=v$ maps $\gamma'$ onto an element of order $m$. As $D$ is isomorphic to a Kleinian group, (\ref{topo giggio}) holds, in this subcase as well, with $Y'=\{1\}$. Thus (\ref{topo giggio}) is established in all cases.

Let $Y$ be the set given by (\ref{topo giggio}), and let us set $Y=J^{-1}(Y')$.
Let $\gamma$ be any element of $P-Y$ and set $\gamma'=J(\gamma)$. Then $\gamma'\in P'-Y'$. Hence for every positive integer $m$ there
exist a Kleinian group $\hatGamma$ and a homomorphism $h':\Gamma'\to \hatGamma$ such that $h'(\gamma)$ has order $m$ in $\hatGamma$. Setting $h= h' \circ J :\Gamma\to\hatGamma$ we find that $h(\gamma)$ has order $m$ in $\hatGamma$. This completes the proof of the first assertion of the proposition.

To prove the second assertion, suppose that $\gamma\in P-Y$ and $m\in\NN$ are given. Fix a Kleinian group $\hatGamma$ and a homomorphism $H:\Gamma\to \hatGamma$ such that $\hatgamma:=H(\gamma)$ has order $m$ in $\hatGamma$. The group $\hatGamma$ is a subgroup of $\zzle(\CC)$, which is in turn isomorphic to a subgroup of ${\rm GL}_2(\CC)$ (Indeed, the adjoint action of $\zzle(\CC)$ on its $3$-dimensional Lie algebra factors through a faithful representation of $\zzle(\CC)$.) In particular $\hatGamma$ is a linear group, and is therefore residually finite according to \cite{malcev}. Since the elements $\hatgamma,\hatgamma^2,\ldots,\hatgamma^{m-1}$ of $\hatGamma$ are nontrivial, there exist a finite group $G$ and a homomorphism $J:\hatGamma\to G$ such that $J(\hatgamma^i)$ is nontrivial for $i=1,\ldots,m-1$. If $h$ denotes the homomorphism $J\circ H:\Gamma\to G$, then $h(\gamma)^i= J(\hatgamma^i)$ is nontrivial for $i=1,\ldots,m-1$, but $h(\gamma)^m$ is trivial since  $\hatgamma^m$ is trivial. Hence $h(\gamma)$ has order $m$ in $G$. This shows that $m$ is a finitistic order for $\gamma$.
\EndProof

\section {Proof of the main theorem}\label{gasoline alley}

This section is devoted to the proof of Theorem \ref{aathur}. 

As in the statement of the theorem, let $\Gamma$ be a finitely generated, torsion-free Kleinian group, let $m>2$ be an integer, and let $X$ denote the set of all elements of $\Gamma$ for which $m$ is a finitistic order. We must show that $\Gamma-X$ is a union of finitely many conjugacy classes. If $\Gamma$ is elementary then it is free abelian by \ref{general stuff reporting for duty}, and by \ref{uneeda} it follows that $X=\Gamma$. For the rest of the proof I will assume that $\Gamma$ is non-elementary.

The required conclusion depends only on the isomorphism class of $\Gamma$. In view of Proposition \ref{no great loss}, we may therefore assume without loss of generality that
\begin{enumerate}
\item $\Gamma$ is geometrically finite,
\item every maximal parabolic subgroup of $\Gamma$ has rank $2$, and 
\item $\Gamma\le\pizzle(K)$ for some number field $K$.
\end{enumerate}

Since $m>1$, the identity element $1$ of $\Gamma$ does not belong to $X$. Since $\Gamma$ is torsion-free and discrete, every non-trivial element of $\Gamma$ is loxodromic or parabolic. Hence we may write $\Gamma-X=Z_\ell\cup Z_p\cup\{1\}$, where $Z_\ell$ (resp. $Z_p$) denotes the set of loxodromic (resp. parabolic) elements of $X$. I shall prove the theorem by showing that each of the sets $Z_\ell$ and $Z_p$ is a union of finitely many conjugacy classes. Note that each of these sets is obviously invariant under conjugation.

Since $\Gamma\le\pizzle(K)$ by (3) above, and since $\Gamma$ is
finitely generated, there is a finite set $S$ of places of $K$ such
that $\Gamma\le\pizzle(\ooks)$. After possibly enlarging $S$ we may
assume that it contains the infinite places, i.e. that it is
admissible. According to \cite[Proposition 3.1.1]{splittings}, there
is a subgroup $\tGamma$ of $\zzle(C)$ such that the quotient
homomorphism $\zzle(\CC)\to\pizzle(\CC)$ restricts to an isomorphism
$Q:\tGamma\to\Gamma$. Since $\Gamma\le\pizzle(\ooks)$, we have
$\tGamma\le \zzle(\ooks)$.

We apply Proposition \ref{congker what}, defining $S$ as above and
using the given value of $m>2$. 
Since  $\tGamma\le \zzle(\ooks)$, Proposition \ref{congker what} implies in particular that there is
a finite set $W\subset\ooks$ such that every element $\gamma$ of $\tGamma$ with $\trace\gamma\notin W$ has $m$ as a finitistic order. Hence $\trace(\tGamma-Q^{-1}(X))\subset W$. 

In particular, if we set $\tZ_\ell=Q^{-1}(Z_\ell)$, then $\trace(\tZ_\ell)$ is contained in $W$ and is therefore finite. On the other hand, since $Z_\ell$ consists of loxodromic elements, we have $2,-2\notin \trace(\tZ_\ell)$. Since $\Gamma$ is geometrically finite by (1) above, it follows from Corollary \ref{cat} that $\tZ_\ell$ contains only finitely many $\Gamma$-conjugacy classes of elements with any given trace. Hence $\tZ_\ell$ contains only finitely many conjugacy classes, and therefore so does $Z_\ell$.

Since $\Gamma$ is geometrically finite, a truncation of $M=\HH^3/\Gamma$ has only finitely many torus components, and hence by Proposition \ref{lame song}, $\Gamma$ has  only finitely many conjugacy classes of maximal parabolic subgroups. Let $k\ge0$ denote the number of these conjugacy classes, and let  $P_1,\ldots,P_k$ be subgroups representing them. By (2) above, each of the $P_i$ has rank $2$. 
According to Proposition \ref{de uddah joik} (and conditions (1) and (2) above), each $P_i$ has a finite subset $Y_i$ such that every positive integer---and in particular $m$---is a finitistic order for every $\gamma\in P_i-Y_i$.  Hence each element of $Z_p$ is conjugate to an element of $\bigcup_{i=1}^kY_i$. This shows that $Z_p$ is contained in a union of finitely many conjugacy classes. 

This completes the proof of Theorem \ref{aathur}.

\bibliographystyle{plain}
\bibliography{finitistic}

\end{document}